\documentclass[12pt]{amsart}

\usepackage{amsmath,amsthm,amssymb}
\font\teneufm=eufm10
\font\seveneufm=eufm7
\font\fiveeufm=eufm5
\newfam\frakturfam

\textfont\frakturfam=\teneufm
\scriptfont\frakturfam=\seveneufm
\scriptscriptfont\frakturfam=\fiveeufm

\newtheorem{df}{Definition}
\newtheorem{lm}{Lemma}
\newtheorem{theor}{Theorem}
\newtheorem{co}{Corollary}

\newtheorem{prob}{Problem}
\def\bee{\begin{eqnarray}}
\def\bes{\begin{eqnarray*}}
\def\eee{\end{eqnarray}}
\def\ees{\end{eqnarray*}}

\def\a{\alpha}
\def\b{\beta}
\def\g{\gamma}

\def\Proof{{\sl Proof.}\ }

\newcommand{\Aut}{\mbox{Aut}}

\title{Universal enveloping algebras and universal derivations of Poisson algebras}

\begin{document}
\date{}
\maketitle

\begin{center}

{\bf Ualbai Umirbaev}\footnote{Supported by an NSF grant DMS-0904713 and
by a grant of Kazakhstan, Department of Mathematics, Eurasian National University,
 Astana, 010008, Kazakhstan and
Department of Mathematics, Wayne State University,
Detroit, MI 48202, USA,
e-mail: {\em umirbaev@math.wayne.edu}}

\end{center}

\begin{abstract}
Let $k$ be an arbitrary field of characteristic $0$.
It is shown that for any $n\geq 1$ the universal enveloping algebras of the  Poisson symplectic algebra $P_n(k)$ and the Weyl algebra $A_n(k)$ are isomorphic and the canonical isomorphism between them easily leads to the Moyal product.
 A basis of the universal enveloping algebra $P^e$ of a free Poisson
algebra $P=k\{x_1,\ldots,x_n\}$ is constructed and proved that the left dependency of a finite number of elements of $P^e$ over $P^e$ is algorithmically recognizable.  We prove that if two elements of a free Poisson algebra do not generate a free two generated subalgebra then they commute. The Fox derivatives on free Poisson algebras are defined and it is proved that an analogue of the Jacobian Conjecture for two generated free Poisson algebras is equivalent to the two-dimensional classical Jacobian Conjecture. A new proof of the tameness of automorphisms of two generated free Poisson algebras is also given.
\end{abstract}

\noindent {\bf Mathematics Subject Classification (2010):} Primary
17B63, 17B40; Secondary 17A36, 16W20.

\noindent

{\bf Key words:} Poisson algebras, universal enveloping algebras, left dependence, automorphisms, derivations.

\section{Introduction}

\hspace*{\parindent}

 The main property of the universal (multiplicative)  enveloping algebra $U_{\mathfrak{M}}(A)$ of an algebra $A$ in a variety of algebras $\mathfrak{M}$ is that the notion of an $A$-bimodule in $\mathfrak{M}$ is equivalent to the notion of a left module over the associative algebra $U_{\mathfrak{M}}(A)$ (see, for example \cite{Jacobson}).
The universal enveloping algebra of a Poisson algebra was studied in \cite{Oh}. A linear basis of the universal enveloping algebra is constructed in \cite{OPC}  for Poisson polynomial algebras, i.e. for Poisson brackets on the polynomial algebras $k[x_1,\ldots,x_n]$. In general case the validity of an analogue of the Poincare-Birkhoff-Witt theorem is still open.

We study the universal enveloping algebras of  Poisson symplectic algebras and free Poisson algebras. For any Poisson algebra $P$ denote by $P^e$ its universal enveloping algebra. Note that the universal enveloping algebra $A^e$ of an associative algebra $A$ in the variety of associative algebras is $A\otimes_k A^{op}$ where $A^{op}$ is an anti-isomorphic copy of $A$ (see, for example \cite{Um11,Um}). It is proved that the universal enveloping algebras of the  Poisson symplectic algebra $P_n$ and the Weyl algebra $A_n$ are isomorphic. The canonical isomorphism between $P_n^e$ and $A_n^e$ naturally leads to the Moyal product (Lemma \ref{l6}).

It is well known that the universal enveloping algebras of free Lie algebras are free associative algebras. P.Cohn  \cite{Cohn} proved that every left ideal of a free associative algebra is a free left module. In fact, it follows from this theorem \cite{Um} that subalgebras of free Lie algebras are free \cite{Shir1,Witt} and automorphisms of finitely generated free Lie algebras are tame \cite{Cohn2}.

Let $P=k\{x_1,\ldots,x_n\}$ be the free Poisson algebra
in the variables $x_1,\ldots,x_n$. We construct a linear basis
 of the universal enveloping algebra $P^e$. Unfortunately left ideals of $P^e$ are not free left $P^e$-modules, i.e., an analogue of Cohn's theorem \cite{Cohn} for free associative algebras is not true in this case. We prove a weaker result which says that the left dependency of a finite set of elements $P^e$ over $P^e$ is algorithmically recognizable. We prove that any two elements of a free Poisson algebra over a field of characteristic $0$ either generate a free Poisson algebra in two variables or commute.  This result was also proved in \cite{MLSh} and an analogue of this result for free associative algebras is well known \cite{Cohn}. We prove that an analogue of the Jacobian Conjecture for two generated free Poisson algebras is equivalent to the two dimensional classical Jacobian Conjecture. A new proof of the tameness of automorphisms of two generated free Poisson algebras is also given.

 This paper is organized as follows. In Section 2 we give the definition of Poisson modules and universal enveloping algebras of Poisson algebras by generators and defining relations.
We define also the universal derivations of Poisson algebras.
Section 3 is devoted to studying the universal enveloping algebras of $P_n$ and $A_n$.
In Section 4 we construct a linear basis
of the universal enveloping algebra $P^e$ of the free Poisson algebra $P=k\{x_1,\ldots,x_n\}$. We introduce a degree function on the universal enveloping algebra $P^e$ and describe its associated graded algebra $gr\,P^e$. In Section 5 we consider the left dependency of elements over $P^e$ and study in details the left dependency of two elements of a special type over $gr\,P^e$. Section 6 is devoted to studying the universal derivations and two generated subalgebras of free Poisson algebras. In Section 7 we give some comments and formulate some open problems.

All vector spaces are considered over an arbitrary fixed field $k$ of characteristic $0$. In the statements of algorithmic character we assume that $k$ is constructive.

\section{Universal enveloping algebras}

\hspace*{\parindent}

Recall that a vector space $P$ over $k$ endowed with two bilinear
operations $x\cdot y$ (a multiplication) and $\{x,y\}$ (a Poisson
bracket) is called {\em a Poisson algebra} if $P$ is a commutative
associative algebra under $x\cdot y$, $P$ is a Lie algebra under
$\{x,y\}$, and $P$ satisfies the following identity (the Leibniz
identity): \bes \{x, y\cdot z\}=\{x,y\}\cdot z + y\cdot \{x,z\}.
\ees

Let $P$ be a Poisson algebra over $k$.
A vector space $V$ over $k$ is called a {\em Poisson module} over $P$
(or {\em Poisson $P$-module}) if there are two bilinear maps
\bes
P\,\times \,V\longrightarrow V, \ \ \ \ ((x,v)\mapsto x\cdot v)
\ees
and
\bes
P\,\times \,V\longrightarrow V, \ \ \ \ ((x,v)\mapsto \{x,v\}),
\ees
such that the relations
\bes
(x\cdot y)\cdot v=x\cdot (y\cdot v),\\
\{\{x,y\},v\}=\{x,\{y,v\}\}-\{y,\{x,v\}\},\\
\{x\cdot y,v\}=y\cdot \{x,v\}+x\cdot \{y,v\},\\
\{x,y\}\cdot v=\{x,y\cdot v\}-y\cdot \{x,v\}
\ees
hold for all $x,y\in P$ and $v\in V$.

The first relation means that every Poisson $P$-module  is a usual (left) module over the  associative and commutative algebra $P$. Sometimes we consider $P$ as an associative and commutative algebra with $x\cdot y$ and a  module over this algebra will be called a module over the commutative algebra $P$.

Let $V$ be a Poisson module over a Poisson algebra $P$.
For every $x\in P$ we denote by $M_x$ the operator of multiplication by $x$ acting on $V$, i.e., $M_x(v)=x\cdot v$ for any $v\in V$.
For every $x\in P$ we define also the Hamiltonian operator  $H_x$ on $V$ by $H_x(v)=\{x,v\}$. Then the Poisson module relations can be written as
\bes
M_{x y}=M_x M_y,\\
  H_{\{x,y\}}=H_x H_y - H_y H_x,\\
  H_{xy}=M_y H_x + M_x H_y, \\
   M_{\{x,y\}}=H_x M_y-M_y H_x,
\ees
respectively.

We consider only Poisson algebras with an identity element  $1$ and unitary modules over Poisson algebras, i.e. $M_1=id$.

Now we give the definition of the universal (multiplicative) enveloping algebra $P^e$ of a Poisson algebra $P$. Let $m_P=\{m_a | a\in P\}$ and $h_P=\{h_a | a\in P\}$ be two copies of the vector space $P$ endowed with two linear isomorphisms $m : P \longrightarrow m_P \ (a\mapsto m_a)$ and $h : P \longrightarrow h_P \ (a\mapsto h_a)$.
Then $P^e$ is an associative algebra over $k$, with an identity $1$, generated by two linear spaces $m_P$ and $h_P$ and defined by the relations
\bes
m_{x y}=m_x m_y, \\
h_{\{x,y\}}=h_x h_y - h_y h_x,\\
h_{xy}=m_y h_x + m_x h_y,\\
 m_{\{x,y\}}=h_x m_y-m_y h_x=[h_x,m_y],\\
 m_1=1
\ees
for all $x,y\in P$.
The last relation comes from our agreement to consider only unitary $P$-modules.
 Note that $h_1=0$.

By the definition of the universal (multiplicative) enveloping algebras \cite{Jacobson} the notion of a bimodule over an algebra is equivalent to the notion of a left module over its universal enveloping algebra.  Let $V$ be an arbitrary Poisson $P$-module. Then $V$ becomes a left $P^e$-module under the action
\bes
m_xv=x\cdot v, h_x v= \{x,v\},
\ees
for all $x\in P$ and $v\in V$. Conversely, if $V$ is a left $P^e$-module then the same formulas turn $V$ to a Poisson $P$-module.
\begin{co}\label{c1}
The category of unitary Poisson modules over a Poisson algebra $P$ and the category of (left) unitary modules over the universal enveloping algebra $P^e$ are equivalent.
\end{co}

The first example of a Poisson $P$-module is $V=P$ under the actions $x\cdot v$ and $\{x,v\}$. Since $1\in P$ and $m_x 1=x$ it follows that the mapping
\bes
m : P\longrightarrow P^e \ \ \ \ \ (x\mapsto m_x)
\ees
is an injection. Therefore we identify $m_x$ with $x$. After this identification the  essential part of the defining relations of $P_e$ are
\bee\label{f1}
h_{\{x,y\}}=h_x h_y - h_y h_x,
\eee
\bee\label{f2}
h_{xy}=y h_x + x h_y,
\eee
\bee\label{f3}
\{x,y\}=h_x y-y h_x=[h_x,y],
\eee
for all $x,y\in P$. From (\ref{f3}) follows that
\bee\label{f4}
[h_x,y]=[x,h_y].
\eee

Let $\Omega_P$ be the left ideal of $P^e$ generated by all $h_x$ where $x\in P$.
Consider the mapping
\bes
\Delta : P \longrightarrow \Omega_P    \ \ \ \ \ (x\mapsto h_x).
\ees
It follows from (\ref{f1}) and (\ref{f2}) that $\Delta$ is a derivation of the Poisson algebra $P$ with coefficients on the Poisson $P$-module $\Omega_P$.

\begin{lm}\label{l1}
$\Delta$ is the universal derivation of the Poisson algebra $P$ and $\Omega_P$ is the universal differential module of $P$.
\end{lm}
\Proof Recall that by the definition it means (see, for example \cite{Um11,Um}) that for any Poisson $P$-module $V$ and for any derivation $d : P \rightarrow V$ there exists a unique homomorphism $\tau : \Omega_P \rightarrow V$ of Poisson $P$-modules such that
$d=\tau \Delta$.

We follow the construction of the universal differential modules by generators and defining relations (see, for example \cite{BergDick}).
Let $\overline{P}$ be a copy of the vector space $P$. Consider the left $P^e$-module
$V=P^e\otimes \overline{P}$. Denote by $\phi : P\rightarrow V$ the mapping defined by $\phi(x)=1 \otimes \overline{x}$ for all $x\in P$. Let $W$ be a submodule of $V$ generated by all
$h_x \otimes\overline{y} - h_y \otimes\overline{x}-1 \otimes\overline{\{x,y\}}$ and $y \otimes \overline{x} + x \otimes\overline{y}-1\otimes\overline{xy}$, where $x,y\in P$. Put $M=V/W$  and denote by $\psi : V \rightarrow M$ the natural homomorphism. It is well known that \cite{BergDick,Um11,Um} the mapping
$D=\psi\phi : P \rightarrow M$ is the universal derivation and $M$ is the universal differential module of $P$.

Since $V$ is a $P^e$-module freely generated by the vector space $\overline{P}$ there exists a $P^e$-module homomorphism $\eta : V\rightarrow \Omega_P$ such that $\eta(\overline{x})=h_x$ for all $x\in P$. It follows from (\ref{f1}) and (\ref{f2}) that $W\subseteq Ker\,\eta$. Consequently, $\eta$ induces the homomorphism
$\theta : M\rightarrow \Omega_P$. Note that from (\ref{f1})--(\ref{f3}) it follows that $\Omega_P$ is a $P^e$-module generated by symbols $h_x$ and defined by relations (\ref{f1}) and (\ref{f2}). Consequently, $\theta$ is an isomorphism. $\Box$

\section{Universal enveloping algebras of $P_n$ and $A_n$}

\hspace*{\parindent}

For an integer $n\geq 1$ the Poisson symplectic algebra $P_n$ is the usual polynomial algebra
$k[x_1,\ldots,x_n,y_1,\ldots,y_n]$ endowed with the Poisson bracket
defined by
\bes
\{x_i,y_j\}=\delta_{ij}, \ \ \{x_i,x_j\}=0, \ \
\{y_i,y_j\}=0,
\ees
where $\delta_{ij}$ is the Kronecker symbol
and $1\leq i,j\leq n$.

By $A_n$ denote the Weyl algebra of index $n$, i.e., $A_n$ is an associative algebra given by generators $X_1,\ldots,X_n,Y_1,\ldots,Y_n$ and defined relations
\bes
[X_i,Y_j]=\delta_{ij}, \ \ [X_i,X_j]=0, \ \
[Y_i,Y_j]=0,
\ees
where $1\leq i,j\leq n$.

Let $k[z_1,z_2,\ldots,z_n]$ be the polynomial algebra in the variables $z_1,z_2,\ldots,z_n$. Denote by
$\mathcal{D}(k[z_1,z_2,\ldots,z_n])$ the algebra of polynomial differential operators on $k[z_1,z_2,\ldots,z_n]$. Note that $A_n$ is isomorphic to $\mathcal{D}(k[z_1,z_2,\ldots,z_n])$:
\bes
Y_i\mapsto z_i, X_i\mapsto \partial_i=\frac{\partial}{\partial z_i}, \ \ \ 1\geq i \geq n.
\ees

The statement of the next lemma was also mentioned in \cite{Oh}.
\begin{lm}\label{l2}
$P_n^e\simeq \mathcal{D}(k[x_1,\ldots,x_n,y_1,\ldots,y_n])\simeq A_{2n}$.
\end{lm}
\Proof It follows directly from (\ref{f1})--(\ref{f4}) that $P_n^e$ is generated by $x_i,y_i,h_{x_i},h_{y_i}$ and satisfies the relations
\bes
[x_i,x_j]=[x_i,y_j]=[y_i,y_j]=0,
[h_{x_i},h_{x_j}]=[h_{x_i},h_{y_j}]\\
=[h_{y_i},h_{y_j}]=0,
[h_{x_i},y_i]={\{x_i,y_i\}}=1, [h_{y_i},x_i]={\{y_i,x_i\}}=-1,
\ees
for all $i,j$.

The Weyl algebra $A_{2n}$ is generated by $X_1,\ldots,X_{2n},Y_1,\ldots,Y_{2n}$ and defined by the relations
\bes
[X_i,Y_j]=\delta_{ij}, \ \ [X_i,X_j]=0, \ \
[Y_i,Y_j]=0,
\ees
where $1\leq i,j\leq 2n$. Consequently, there exists a surjective homomorphism $\phi : A_{2n}\longrightarrow P_n^e$ such that
\bes
X_{i}\mapsto x_i,X_{i+n}\mapsto y_i,Y_{i}\mapsto h_{y_i}, Y_{i+n}\mapsto -h_{x_i}
\ees
for all $1\leq i\leq n$. Since $A_{2n}$ is simple it follows that $\phi$ is an isomorphism. $\Box$

Now let $A$ be an arbitrary associative algebra over $k$ with an identity. Then the universal enveloping algebra $A^e$ of $A$ is $A\otimes_k A^{op}$, where $A^{op}$ is an algebra anti-isomorphic to $A$. Algebra $A^{op}$ is considered together with an anti-isomorphism $' : A\longrightarrow A^{op}$. In fact, each  operator of the left multiplication $l_x$ is represented by $x\otimes 1$ and each operator of the right multiplication $r_x$ is represented by $1\otimes x'$ in $A\otimes_k A^{op}$. Following the main property of the universal enveloping algebras (see \cite{Jacobson}) we get that the category of associative unitary $A$-bimodules and the category of left unitary $A^e$-modules are equivalent.

Note that $A_n^{op}$ is isomorphic to $A_n$ again. An isomorphism, can be chosen, for example, as $\varphi : A_n\longrightarrow A_n^{op}$ with $\varphi(X_i)=X_i', \varphi(Y_i)=-Y_i'$ for all $i$. Then,
$A_n^e=A_n\otimes_k A_n^{op}\simeq A_n\otimes_k A_n\simeq A_{2n}$ and consequently, $P_n^e\simeq A_n^e$. We proved

\begin{theor}\label{t1} For every integer $n\geq 1$ the universal enveloping algebras of the Poisson symplectic algebra $P_n$ and the Weyl algebra $A_n$ are isomorphic to the Weyl algebra $A_{2n}$, i.e.
\bes
P_n^e\simeq A_n^e\simeq A_{2n}.
\ees
\end{theor}
\begin{co}\label{c2}
The category of unitary Poisson modules over the Poisson symplectic algebra $P_n$ and the category of unitary bimodules over the Weyl algebra $A_n$ are equivalent.
\end{co}

The notion of bimodules in the case of Poisson algebras corresponds to the notion of Poisson modules, as it follows from the commutativity and anti-commutativity of operations.

Therefore $P_n^e$ is isomorphic to $A_n^e$. We use this isomorphism to get the Moyal product. In fact, there are many isomorphisms between $A_n^e$ and $P_n^e$. We choose "the most canonical" one between them in the next lemma.
\begin{lm}\label{l3} There exists a unique isomorphism  $\theta : A_n^e\longrightarrow P_n^e$ such that
\bes
\theta(X_i\otimes 1)=x_i+1/2h_{x_i}, \theta(1\otimes X_i')=x_i-1/2h_{x_i},\\
\theta(Y_i\otimes 1)=y_i+1/2h_{y_i}, \theta(1\otimes Y_i')=y_i-1/2h_{y_i},
\ees
for all $1\leq i\leq n$.
\end{lm}
\Proof The existence of $\theta$ follows from (\ref{f1})--(\ref{f3}). For example,
\bes
[\theta(X_i\otimes 1),\theta(Y_i\otimes 1)]=[x_i+1/2h_{x_i},y_i+1/2h_{y_i}]\\
=1/2[x_i,h_{y_i}]+1/2[h_{x_i},y_i]+1/4[h_{x_i},h_{y_i}]= \{x_i,y_i\}=1, \\ \
  [\theta(1 \otimes X_i'),\theta(1 \otimes Y_i')]=[x_i-1/2h_{x_i},y_i-1/2h_{y_i}]\\ \
=-1/2[x_i,h_{y_i}]-1/2[h_{x_i},y_i]+1/4[h_{x_i},h_{y_i}]= -\{x_i,y_i\}=-1, \\ \
  [\theta(X_i\otimes 1),\theta(X_j\otimes 1)]=[x_i+1/2h_{x_i},x_j+1/2h_{x_j}]\\
=1/2[x_i,h_{x_j}]+1/2[h_{x_i},x_j]+1/4[h_{x_i},h_{x_j}]= 0, \\
\ees
if $i\neq j$. Obviously, $\theta$ is surjective. And it is also injective since $A_n^e$ is simple.  $\Box$

Denote by $L$ the linear space $kx_1+\ldots+kx_n+ky_1+\ldots+ky_n$.
We define a linear bijection
\bes
w :  P_n \rightarrow A_n
\ees
by
\bes
w(l_1l_2\ldots l_p)= 1/p! \sum_{\pi\in S_p} l_{\pi(1)}l_{\pi(2)}\ldots l_{\pi(p)}
\ees
for all $l_1,l_2,\ldots,l_p\in kx_1+\ldots+kx_n+ky_1+\ldots+ky_n$. This mapping is called the symmetrization \cite{Dixmier}

Let us introduce some notations. Let $\mathbb{Z}_{+}$ be the set of all nonnegative integers. For every $\a = (i_1,i_2,\ldots,i_t)\in \mathbb{Z}_{+}^t$ we put $|\a|=i_1+i_2+\ldots+i_t$ and $\a!=i_1!i_2!\ldots i_t!$.
If $b=(b_1,b_2,\ldots,b_t)$, where $b_1,b_2,\ldots,b_t$ is a subset of commuting elements of an algebra $B$, then we put
\bes
b^{\a}=b_1^{i_1}b_2^{i_2}\ldots b_t^{i_t}.
\ees
We set
\bes
\partial_x=(\frac{\partial}{\partial x_1}, \frac{\partial}{\partial x_2},\ldots,\frac{\partial}{\partial x_n}),
\partial_y=(\frac{\partial}{\partial y_1},\frac{\partial}{\partial y_2},\ldots, \frac{\partial}{\partial y_n}),\\
\partial=(\frac{\partial}{\partial x_1},\ldots,\frac{\partial}{\partial x_n},
\frac{\partial}{\partial y_1},\ldots,\frac{\partial}{\partial y_n}), \\
h_x=(h_{x_1},h_{x_2},\ldots,h_{x_n}), h_y=(h_{y_1},h_{y_2},\ldots,h_{y_n}),\\
h=(h_{x_1},\ldots,h_{x_n},h_{y_1},\ldots,h_{y_n}),
\ees
where $\frac{\partial}{\partial x_i}, \frac{\partial}{\partial y_i}$ are usual partial derivatives of $P_n$.

\begin{lm}\label{l4} If $f\in P_n$ then
\bes
\theta(w(f)\otimes 1)=\sum_{\g\in \mathbb{Z}_{+}^{2n}} \frac{1}{\g!2^{|\g|}}\partial^{\g}(f)h^{\g}
= \sum_{\g\in \mathbb{Z}_{+}^{2n}} \frac{1}{\g!2^{|\g|}}\partial^{\g}(f)h^{\g}.
\ees
\end{lm}
\Proof
By the linearity, it is sufficient to prove the statement of the lemma only for elements of the form $f=l^p$, where $l\in kx_1+\ldots+kx_n+ky_1+\ldots+ky_n$. In this case $w(f)=w(l)^p$. Hence
\bes
\theta(w(f)\otimes 1)=\theta(w(l)^p\otimes 1)=(\theta(w(l)\otimes 1))^p
=(l+1/2h_l)^p.
\ees
Note that $l$ and $h(l)$ commute. Consequently,
\bes
\theta(w(f)\otimes 1)=\sum_{i=0}^{p}
\left(\begin{array}{cc}
p\\
i\\
\end{array}\right)
 \frac{1}{2^i} l^{p-i}h_l^i=\sum_{i=0}^{p}
 \frac{1}{i!2^i} \frac{\partial^i f}{\partial l^i}h_l^i.
\ees
Let $l=a_1x_1+\ldots+a_nx_n+b_1y_1+\ldots+b_ny_n$ and put $c=(a_1,\ldots,a_n,b_1,\ldots,b_n)$. The elements $h_{x_1},\ldots,h_{x_n},h_{y_1},\ldots,h_{x_n}$ commute. Using multinomial coefficients we can write
\bes
h_l^i=\sum_{\lambda\in \mathbb{Z}_{+}^{2n}, |\lambda|=i} \frac{i!}{\lambda !} c^{\lambda} h^{\lambda}.
\ees
Note that $\partial^{\lambda}(f)=c^{\lambda}\frac{\partial^i f}{\partial l^i}$ if $|\lambda|=i$. Consequently,
\bes
\theta(w(f)\otimes 1)=\sum_{\g\in \mathbb{Z}_{+}^{2n}} \frac{1}{\g!2^{|\g|}}\partial^{\g}(f)h^{\g}.
\ees
The second equality of the lemma can be proved similarly. $\Box$

If  $\a=(i_1,\ldots,i_n,j_1,\ldots,j_n)\in \mathbb{Z}_{+}^{2n}$ then we put
$\a_1=(i_1,\ldots,i_n)\in \mathbb{Z}_{+}^{n}$ and $\a_2=(j_1,\ldots,j_n)\in \mathbb{Z}_{+}^{n}$. We often write $\a=(\a_1,\a_2)$ and put $\a^*=(\a_2,\a_1)$.

\begin{lm}\label{l5} If $f\in P_n$ then
\bes
\frac{1}{\g!}h^{\g}f=\sum_{\a+\b=\g}
\frac{(-1)^{|\a_2|}}{\a!\b!} \partial^{\a^*}(f) h^{\b}
\ees
for all $\g\in \mathbb{Z}_{+}^{2n}$.
\end{lm}
\Proof  Put $H_{a}(f)=[h_a,f]$ for all $a, f\in P$ and put also $H=(h_{x_1},\ldots,h_{x_n},h_{y_1},\ldots,h_{y_n})$. Using the relations
\bes
h_{x_i}f=[h_{x_i},f]+fh_{x_i}, h_{y_i}f=[h_{y_i},f]+fh_{y_i},
\ees
and multinomial coefficients, we get
\bes
\frac{1}{\g!}h^{\g}f=\sum_{\a+\b=\g}
\frac{1}{\a!\b!} H^{\a}(f) h^{\b}.
\ees
Note that
\bes
H_{x_i}(f)=[h_{x_i},f]=\frac{\partial}{\partial y_i}(f), H_{y_i}(f)=[h_{y_i},f]=-\frac{\partial}{\partial x_i}(f)
\ees
for all $i$. Consequently, $H^{\a}(f)=(-1)^{|\a_2|} \partial^{\a^*}(f)$.
 $\Box$

Consider the mapping
\bes
\rho_w : P_n\longrightarrow P_n^e, \ \ \  f\mapsto w(f)\otimes 1, \ \ f\in P_n.
\ees
\begin{lm}\label{l6} The product
\bes
f*_wg=\rho_w^{-1}(\rho_w(f)\rho_w(g))
\ees
is the Moyal product on $P_n$, i.e.
\bes
f*_wg=
\sum_{\a\in \mathbb{Z}_{+}^{2n}} \frac{(-1)^{|\a_2|}}{\a!2^{|\a}} \partial^{\a}(f)\partial^{\a^*}
(g).
\ees
\end{lm}
\Proof Direct calculation gives
\bes
\theta(w(f)\otimes 1)\theta(w(g)\otimes 1)
=\sum_{\g,\delta\in \mathbb{Z}_{+}^{2n}} \frac{1}{\g!\delta!2^{|\g+\delta|}}\partial^{\g}(f)h^{\g}
\partial^{\delta}(g)h^{\delta}\\
=\sum_{\g,\delta\in \mathbb{Z}_{+}^{2n}} \sum_{\a+\b=\g}
\frac{(-1)^{|\a_2|}}{\a!\b!\delta!2^{|\g+\delta|}}\partial^{\g}(f)\partial^{\a^*+\delta}
(g)h^{\b+\delta}\\
=\sum_{\mu}(\sum_{\b+\delta=\mu} \sum_{\a} \frac{(-1)^{|\a_2|}}{\a!\b!\delta!2^{|\a+\b+\delta|}} \partial^{\a+\b}(f)\partial^{\a^*+\delta}
(g))h^\mu\\
= \sum_{\mu}
\frac{1}{\mu!2^{|\mu|}} \partial^{\mu}
(\sum_{\a} \frac{(-1)^{|\a_2|}}{\a!2^{|\a}} \partial^{\a}(f)\partial^{\a^*}
(g))h^\mu=
\sum_{\mu}
\frac{1}{\mu!2^{|\mu|}} \partial^{\mu}
(f*_wg)h^\mu.
\ees
Consequently,
\bes
\rho_w^{-1}(\rho_w(f)\rho_w(g))=\rho_w^{-1}(\theta(w(f)\otimes 1)\theta(w(g)\otimes 1))=
f*_wg
\ees
for all $f,g\in P_n$.  $\Box$

\section{Enveloping algebras of free Poisson algebras}

\hspace*{\parindent}

Let $g$ be a Lie algebra with a linear
basis $e_1,e_2,\ldots,e_k,\ldots$. The Poisson symmetric algebra $PS(g)$ of
$g$ is the usual polynomial algebra $k[e_1,
e_2,\ldots,e_k,\ldots]$ endowed with the Poisson bracket defined
by  \bes \{e_i,e_j\}=[e_i,e_j] \ees for all $i,j$, where $[x,y]$
is the multiplication in the Lie algebra $g$.

Denote by $P=k\{x_1,x_2,\ldots,x_n\}$ the free Poisson algebra over $k$ in the variables $x_1,x_2,\ldots,x_n$.
From now on let $g=Lie\langle x_1,x_2,\ldots,x_n\rangle$ be the free Lie algebra with free (Lie)
generators $x_1,x_2,\ldots,x_n$. It is well-known (see, for
example \cite{Shest}) that the Poisson symmetric algebra $PS(g)$ is the free Poisson algebra $P=k\{x_1,x_2,\ldots,x_n\}$ in the variables $x_1,x_2,\ldots,x_n$.

By $\deg$ we denote the
standard degree function of the homogeneous algebra $P$, i.e.
$\deg(x_i)=1$, where $1\leq i \leq n$. Note that
\begin{eqnarray*}
\deg\,\{f,g\}= \deg\,f+\deg\,g
\end{eqnarray*}
if $f$ and $g$ are homogeneous and $\{f,g\}\neq 0$. By
$\deg_{x_i}$ we denote the degree function on $P$ with respect to
$x_i$.
If $f$ is homogeneous with respect to each $\deg_{x_i}$, where
$1\leq i \leq n$, then $f$ is called multihomogeneous.

Let us choose a multihomogeneous linear basis \bes
x_1,x_2,\ldots,x_n,\,[x_1,x_2],\ldots,[x_1,x_n],\ldots,[x_{n-1},x_n],\,[[x_1,x_2],x_3],\ldots
\ees of a free Lie algebra $g$ and denote the elements of this
basis by
\bee
\label{f5} e_1, e_2, \ldots, e_m, \ldots.
\eee

The algebra $P=k\{ x_1,x_2,\ldots,x_n\}$ coincides with
the polynomial algebra on the elements (\ref{f5}).
Consequently, the set of all words of the form
\bee\label{f6}
u=e^{\a}=e_1^{i_1}e_2^{i_2}\ldots e_m^{i_m},
\eee
where $0\leq i_k$, $1\leq k\leq m$, and $m\geq 0$, forms a linear basis of $P$.
The basis (\ref{f6}) is multihomogeneous since so is (\ref{f5}).

Let $k\{x_1,x_2,\ldots,x_n,y\}$ be the free Poisson algebra in the variables $x_1,x_2,\ldots,x_n,y$. Denote by $W$ the set of all homogeneous of degree one with respect to $y$ elements of
$k\{x_1,x_2,\ldots,x_n,y\}$.

\begin{theor}\label{t2}
Let $P=k\{x_1,x_2,\ldots,x_n\}$ be the free Poisson algebra over a field $k$ in the variables $x_1,x_2,\ldots,x_n$ and let $P^e$ be its universal enveloping algebra. Then the following statements are true:

(i) The subalgebra $A$ (with identity) of $P^e$ generated by $h_{x_1},h_{x_2},\ldots, h_{x_n}$ is the free associative algebra in the variables $h_{x_1},h_{x_2},\ldots, h_{x_n}$;

(ii) The left commutative $P$-module $P^e$ is isomorphic to the left commutative $P$-module $P\otimes_k A$.
\end{theor}
\Proof Recall that $P\{x_1,x_2,\ldots,x_n,y\}$ is the Poisson symmetric algebra of the free Lie algebra $Lie\langle x_1,x_2,\ldots,x_n,y\rangle$. The elements  of the form
\bes
\{x_{i_1},\{x_{i_2},\ldots,\{x_{i_k},y\}\ldots\}\}= h_{x_{i_1}}h_{x_{i_2}}\ldots h_{x_{i_k}}(y)
\ees
are linearly independent in the free Lie algebra $Lie\langle x_1,x_2,\ldots,x_n,y\rangle$.
Consequently, the elements of the form
\bee\label{f7}
h_{x_{i_1}}h_{x_{i_2}}\ldots h_{x_{i_k}}
\eee
 are linearly independent in $P^e$.

  Using (\ref{f1})--(\ref{f4}), it can be easily shown that every element of $P^e$ can be written as a linear combination of elements $p w$, where $p\in P$ and $w$ is an element of the form (\ref{f7}).

Let $B_1$ be the linear basis (\ref{f5})  of the free Lie algebra $g=Lie\langle x_1,x_2,\ldots,x_n\rangle$. Denote by $B_2$ the set of all elements of the form $wy$, where $w$ is an element of the form (\ref{f7}). Note that the set of elements $B_1\bigcup B_2$ is linearly independent. We can choose a set of elements $B_3$ of degree $\geq 2$ in $y$ such that $B_1\bigcup B_2\bigcup B_3$ is a linear basis of $Lie\langle x_1,x_2,\ldots,x_n,y\rangle$.
Then $P\{x_1,x_2,\ldots,x_n,y\}$ is a polynomial algebra in the set of variables $B_1\bigcup B_2\bigcup B_3$. Consequently, $W$ is a free left module over the polynomial algebra $k[B_1]$ and $B_2$ is a set of free generators. Note that $P=k[B_1]$. $\Box$

\begin{co}\label{c3}
Every nonzero element $u$ of the universal enveloping algebra $P^e$ can be uniquely written in the form
\bee\label{f8}
u=\sum_{i=1}^k p_i w_i,
\eee
where $0\neq p_i\in P$ for all $i$ and $w_1,w_2,\ldots, w_k$ are different elements of the form (\ref{f7}).
\end{co}

Put $h_{x_i}<h_{x_j}$ if $i<j$. Let $u,v$ be two elements of the form (\ref{f7}). Then put $u<v$ if $\deg\,u<\deg\,v$ or $\deg\,u=\deg\,v$ and $u$ precedes $v$ in the lexicographical order.  We can assume $w_1<w_2<\ldots <w_k$ in (\ref{f8}). Then $w_k$ is called the {\em leading monomial} of $u$ and $p_k$ is called the {\em leading coefficient} of $u$. We will write $w_k=ldm(u)$ and $p_k=ldc(u)$. The {\em leading term} of $u$ is defined by $ldt(u)=ldc(u) ldm(u)$.

\begin{lm}\label{l7}
If $u$ and $v$ are arbitrary nonzero elements of $P^e$ then $ldc(uv)=ldc(u)ldc(v)$ and $ldm(uv)=ldm(u)ldm(v)$.
\end{lm}
\Proof
Note that if $u$ and $v$ are two elements of the form (\ref{f8}) then to put the product $uv$ into the form (\ref{f8}) again we need to use only the relations (\ref{f3}). This means that $h_{x_i}$ and $y\in P$ commutes modulo terms of smaller degrees in the variables $h_{x_1},h_{x_2},\ldots,h_{x_n}$. Consequently, we can put $uv$ into the form (\ref{f8}) with the leading monomial $ldm(u)ldm(v)$ and the leading coefficient $ldc(u)ldc(v)$.
 $\Box$

Now we introduce a degree function $hdeg$ (or $h$-degree function) on $P^e$.
Let $u$ be an element of $P^e$ written in the form (\ref{f8}). Then we put $hdeg\,u= max_{i=1}^k deg\,w_i$ and $hdeg\,0=-\infty$.  We say that $u$ is homogeneous with respect to $hdeg$ if $deg\,w_1=deg\,w_2=\ldots =deg\,w_k$.
It follows directly from Lemma \ref{l6}  and (\ref{f3}) that
\bes
hdeg\,uv=hdeg\,u+hdeg\,v
\ees
for every $u$ and $v$ from $P^e$, i.e., $hdeg$ is a degree function on $P^e$.  Denote by $\overline{u}$ the highest homogeneous part of $u$ with respect to $hdeg$.

Denote by $U_i$ the subset of all elements $u$ of $P^e$ with $hdeg\,u\leq i$. Then,
\bes
P=U_0\subset U_1\subset U_2\subset\ldots \subset U_k\subset \ldots ,
\ees
is a filtration of $P^e$, i.e., $U_i U_j\subseteq U_{i+j}$ for all $i,j\geq 0$. Put also
\bes
gr\,P^e= gr\,U_0\oplus gr\,U_1\oplus gr\,U_2\oplus\ldots \oplus gr\,U_k\oplus\ldots,
\ees
where $gr\,U_0=P$ and $gr\,U_i=U_i/U_{i-1}$ for all $i\geq 1$.
Denote by $\varphi_i : U_i\rightarrow gr\,U_i$ the natural projection for every $i\geq 1$ and put $\varphi_0=id : P\rightarrow P$. We define also
\bee\label{f9}
\varphi=\{\varphi_i\}_{i\geq 0} : P^e\rightarrow gr\,P^e
\eee
 by $\varphi(u)=\varphi_i(u)$ if $u\in U_i\setminus U_{i-1}$ for every $i\geq 1$ and $\varphi(u)=u$ if $u\in P$.

The multiplication of $P^e$ induces a multiplication on $gr\,P^e$ and the graded vector space $gr\,P^e$ becomes an algebra.

Consider $B=P\otimes_k A$ as a tensor product of associative algebras. Then $B$ is a free associative algebra over $P$ in the variables $h_{x_1},h_{x_2},\ldots,h_{x_n}$.

\begin{theor}\label{t3}
The graded algebra $gr\,P^e$ is isomorphic to $B=P\otimes_k A$.
\end{theor}
\Proof By (\ref{f3}), $P$ is in the center of the algebra $gr\,P^e$ and $gr\,P^e$ is generated by $\varphi(h_{x_1}),\varphi(h_{x_2}),\ldots,\varphi(h_{x_n})$ as an algebra over $P$. Note that $B=P\otimes_k A$ is a free associative algebra over $P$. Then there is a $P$-algebra homomorphism $\psi : B \rightarrow gr\,P^e$ such that $\psi(h_{x_i})=\varphi(h_{x_i})$ for all $i$.

Let $T_s$ be the space of homogeneous with respect to $hdeg$ elements of $P^e$ of degree $s\geq 1$ and $B_s$ be the space of homogeneous of degree $s$ elements of $B$. There is an obvious isomorphism between the spaces $T_s$ and $A_s$ established by the $P$-module homomorphism in Theorem \ref{t2}. Note that $U_s=U_{s-1}+T_s$, $U_s/U_{s-1}\simeq T_s\simeq B_s$, and $\psi_{|B_s} : B_s\rightarrow gr\,U_s$ is an isomorphism of $P$-modules.  Consequently, $Ker\,(\psi)=0$ and $\psi$ is an isomorphism of algebras. $\Box$

\section{The left dependency}

\hspace*{\parindent}

We use the notations of Section 4.
\begin{lm}\label{l8} Let $u\in P^e$ and $hdeg\,u=m$. Then there exists $v\in P^e$ such that $\lambda^{m+1} u= v\lambda$.
\end{lm}
\Proof By (\ref{f3}), $\lambda u=u \lambda +u_1$, where $u_1$ has degree less than $m$. Consequently, $\lambda^{m+1} u = w \lambda$ by induction on $m$. $\Box$

Let $u$ and $v$ be two elements of the form (\ref{f7}).  We write $u\ll v$ if u is a left divisor of $v$, i.e., $v=tu$ for some $t$ of the form (\ref{f7}).

\begin{df}\label{d1} Let $u$ and $v$ be nonzero elements of $P^e$ such that $ldm(u)=ldm(v)$. Let $r=gcd(ldc(u),ldc(v))$ be the greatest common divisor of $ldc(u)$ and $ldc(v)$. Then
put
\bes
(u,v)_c= (ldc(v)/r) u - (ldc(u)/r) v.
\ees
\end{df}
Note that $(u,v)_c=0$ or $ldm((u,v)_c)<ldm(u)=ldm(v)$.

\begin{lm}\label{l9}
Let $s_1,s_2,\ldots,s_k$ be a finite set of nonzero elements of $P^e$. If $ldm(s_i)$ and $ldm(s_j)$ are not comparable with respect to $\ll$ for every $i\neq j$ then the elements $s_1,s_2,\ldots,s_k$ are left independent over $P^e$.
 \end{lm}
\Proof
Suppose that
\bee\label{f10}
\sum_{r=1}^k u_r s_r=0.
\eee
By Lemma \ref{l7}, $ldm(u_r s_r)=ldm(u_r) ldm(s_r)$ for every $r$. Suppose that the equation (\ref{f10}) is not trivial, i.e., at least one of the coefficients $u_r$ is nonzero. Then,  comparing the leading monomials of the summands, we conclude that $ldm(u_i) ldm(s_i)=ldm(u_j) ldm(s_j)\neq 0$ for some $i\neq j$. It is possible if and only if $ldm(s_i)\ll ldm(s_j)$ or $ldm(s_i)\gg ldm(s_j)$. This contradicts to the condition of the lemma. $\Box$

\begin{lm}\label{l10}
 Let $s_1,s_2,\ldots,s_k$ be a finite set of nonzero elements of $P^e$. Suppose that $ldm(s_i)\gg ldm(s_j)$ and $ldm(s_i)= t\cdot ldm(s_j)$. Put $s_i'=(s_i,s_j)_c$.  Then the elements $s_1,s_2,\ldots,s_k$  are left dependent over $P^e$ if and only if the elements \\ $s_1,s_2,\ldots,s_{i-1},s_i',s_{i+1},\ldots,s_k$ are left dependent over $P^e$.
\end{lm}
\Proof Consider the equation (\ref{f10}).
By Lemma \ref{l8}, there exists a number $m$ such that $\lambda^m u_i=v_i \lambda$ for every $\lambda\in P$ and $i$.
We choose
\bes
\lambda=ldc(s_j)/\mu_{ij},
\ees
where $\mu_{ij}=gcd(ldc(s_i),ldc(s_j))$.

Note that
\bes
\lambda s_i=s_i'+(ldc(s_i)/\mu_{ij}) s_j.
\ees
Then (\ref{f10}) is equivalent to
\bes
\sum_{r=1}^k v_r\lambda s_r=\sum_{r\neq i,j} v_r \lambda s_r+ v_i s_i'+(v_j\lambda + v_i (ldc(s_i)/\mu_{ij})) s_j=0.
\ees
Obviously, this relation is trivial if and only if (\ref{f10}) is trivial.
$\Box$
\begin{theor}\label{t4} In the universal enveloping algebra $P^e$ of the free Poisson algebra
$P=k\{x_1,x_2,\ldots,x_n\}$
the left dependency of a finite system of elements is algorithmically recognizable.
\end{theor}
\Proof Let $s_1,s_2,\ldots,s_k$ be a finite set of nonzero elements of $P^e$. If $ldm(s_i)$ and $ldm(s_j)$ are not comparable with respect to $\ll$ for every $i\neq j$ then the elements $s_1,s_2,\ldots,s_k$ are left independent over $P^e$ by Lemma \ref{l9}.

If $ldm(s_i)\gg ldm(s_j)$ for some $i\neq j$ then we can change $s_i$ by $s_i'$, according to Lemma \ref{l10}. Note that $ldm(s_i)>ldm(s_i')$. If $s_i'=0$ then the elements $s_1,s_2,\ldots,s_k$ are left dependent.  If $s_i'\neq 0$ then we will apply the same discussions to the new system of elements $s_1,s_2,\ldots,s_{i-1},s_i',s_{i+1},\ldots,s_k$.  This process stabilizes after a finite number of steps since the set of all leading monomials is well-ordered. $\Box$

 Recall that $B=P\otimes_k A$ is the free associative algebra over $P$ in the variables $y_1=h_{x_1}, y_2=h_{x_2},\ldots,y_n=h_{x_n}$. Considering $B$ as an algebra over $P$, we define a degree function $Deg$ on $B$ by $Deg\,y_i=1$, where  $1\leq i\leq n$.

Denote by $L$ the Lie subalgebra of the Lie algebra $A^{(-)}$ generated by $y_1,y_2,\ldots,y_n$. Then $L$ is  a free Lie algebra and $y_1,y_2,\ldots,y_n$ are free
generators of $L$. We need the next purely Lie-theoretical statement.
\begin{lm}\label{l11}
Let $f$ be a homogeneous nonlinear element of $L$ such that $f=f_1y_1+f_2y_2+\ldots+f_ny_n$ and $f_n\neq 0$ in $A$.
Then there exists $i\leq n-1$ such that $ldm(f_i)> ldm(f_n)$.
\end{lm}
\Proof
Recall that a nonempty associative word $u$ in the alphabet $y_1,y_2,\ldots,y_n$ is called a {\em Lindon-Shirshov} word (see, for example \cite{Bokut}) if for every nonempty words $v$ and $w$ the equality $u=vw$ implies $vw>wv$. It is well known that if $f\in L$ then $ldm(f)$ is a Lindon-Shirshov word \cite{Bokut}.

Suppose that for every nonzero $f_i$, where $i\leq n-1$, $ldm(f_i)\leq ldm(f_n)$, then  $ldm(f)=ldm(f_n)y_n=u$  and $u$ is a Lindon-Shirshov word. Put $v=ldm(f_n)$, then  $vy_n>y_nv$. Since $y_n$ is the greatest symbol of the alphabet it follows that $v=y_nw$ and $wy_n>y_nw$. Continuing the same discussions we can get that $u=y_n^s$ for some $s\geq 2$.
Note that $y_n^s$ is not a Lindon-Shirshov word. $\Box$
\begin{lm}\label{l12} Let $f$ and $g$ be nonzero homogeneous with respect to $Deg$ elements of\\
 $P\otimes_k L\subseteq P\otimes_k A=B$. If $f$ and $g$ are left dependent over $B$ then $Deg\,f=Deg\,g$.
\end{lm}
\Proof
Suppose that $Deg\,f\geq Deg\,g$. The elements $f$ and $g$ are left dependent in a free associative algebra over $P$. Then there exist a nonzero $\lambda\in P$ and a homogeneous of degree $Deg\,f-Deg\,g$ element $T\in B$ such that $\lambda f=T g$. Changing $f$ by $\lambda f$ we may assume that
\bee\label{f11}
f=T g.
\eee
Every nonzero homogeneous element $b\in B$ can be represented as
\bes
b= \b_1\otimes a_1+\b_2\otimes a_2+\ldots +\b_s\otimes a_s
\ees
with the least possible $s$, $ldm(a_1)<ldm(a_2)<\ldots <ldm(a_s)$, and $ldc(a_1)=ldc(a_2)=\ldots =ldc(a_s)=1$. We call this representation of $b$ a {\em short representation}. From the minimality of $s$ it follows that $\b_1,\b_2,\ldots,\b_s$ are linearly independent elements of $P$ and
$a_1,a_2,\ldots,a_s$ are linearly independent elements of $A$. If $b\in P\otimes_k L$ then we can easily find a short representation of $b$ with an additional condition $a_i\in L$ for all $i$.
Let
\bes
f=\a_1\otimes l_1+\a_2\otimes l_2+\ldots+\a_s\otimes l_s, \\
g=\b_1\otimes m_1+\b_2\otimes m_2+\ldots+\b_t\otimes m_t, \\
T=\g_1\otimes W_1+\g_2\otimes W_2+\ldots+\g_r\otimes W_r
\ees
be short representations of $f$, $g$, $T$ such that $l_i, m_j\in L$ for all $i,j$. Then,
\bes
f=\sum_{i,j} \g_i\b_j W_i m_j.
\ees
From this representation of $f$ we can get every other short representation of $f$ by linear transformations over $k$.
Consequently, $l_i$ belongs to the left ideal of $A$ generated by $m_1,m_2,\ldots,m_t$ for all $i$.
It follows from (\ref{f11}) that $\a_s= \g_r \b_t$ and $ldm(l_s)=ldm(W_r) ldm(m_t)$. Consequently,
\bes
l_s=g_1m_1+g_2m_2+\ldots+g_tm_t
\ees
 and $ldm(g_1),ldm(g_2),\ldots,ldm(g_{t-1})\leq ldm(g_t)=W_r$.

It follows from \cite{Um6,Um9} that $l_s$ belongs to the Lie subalgebra generated by $m_1,m_2,\ldots,m_t$. We fix $h\in Lie<z_1,z_2,\ldots,z_t>$ such that $l_s=h(m_1,m_2,\ldots,m_t)$. We may assume that $h$ is homogeneous since
$m_1,m_2,\ldots,m_t$ are homogeneous and have the same degrees. In this case their linear independency implies \cite{Um} their freeness in the Lie algebra $L$. Hence the left $A$-submodule generated by these elements is free with the same free generators \cite{Um}. If $h=h_1z_1+h_2z_2+\ldots+h_tz_t$ then
$h_i(m_1,m_2,\ldots,m_t)=g_i$ for all $i$. Put $z_1<z_2<\ldots<z_t$. Note that
\bes
ldm(h_i(m_1,m_2,\ldots,m_t))=ldm(h_i)(ldm(m_1),ldm(m_2),\ldots,ldm(m_t))
\ees
since $ldm(m_1)<ldm(m_2)<\ldots <ldm(m_t)$ and have one and the same degrees. The same idea gives $ldm(h_i)<ldm(h_j)$  if and only if $ldm(g_i)<ldm(g_j)$. Then,
\bes
ldm(h_1),ldm(h_2),\ldots,ldm(h_{t-1})\leq ldm(h_t).
\ees
If $h$ is not linear then this contradicts to Lemma \ref{l11}, therefore $h$ is linear. It is possible if and only if $Deg\,T=0$, hence $Deg\,f=Deg\,g$. $\Box$

\section{Universal derivations}

\hspace*{\parindent}

As before, $P=k\{x_1,x_2,\ldots,x_n\}$ is the free Poisson algebra in the variables $x_1,x_2,\ldots,x_n$. Denote by $\Omega_P$ the left ideal of $P^e$ generated by $h_{x_1},h_{x_2},\ldots, h_{x_n}$. By Theorem \ref{t2},
\bes
\Omega_P = P^eh_{x_1}\oplus P^eh_{x_2}\oplus\ldots \oplus P^eh_{x_n},
\ees
i.e., $\Omega_P$ is a free left $P^e$-module. Note that
\bes
P^e=P \oplus \Omega_P.
\ees

Consider
\bes
H : P \rightarrow \Omega_P
\ees
such that $H(p)=h_p$ for all $p\in P$. By Lemma \ref{l1}, $H$ is the universal derivation of $P$ and $\Omega_P$ its universal differential module.

Recall that a set of elements $f_1,f_2,\ldots,f_k$ of the free Poisson algebra $P$ is called {\em Poisson free} or {\em Poisson independent} if the Poisson subalgebra of $P$ generated by these elements is the  free Poisson algebra with free generators $f_1,f_2,\ldots,f_k$. Otherwise these elements are called {\em Poisson dependent}.

\begin{lm}\label{l13}
 Let $f_1,f_2,\ldots,f_k$ be arbitrary elements of the free Poisson algebra
 $P$ over a field $k$ of characteristic $0$.
If the elements $f_1,f_2,\ldots,f_k$ are Poisson dependent then the elements
$H(f_1),H(f_2),\ldots,H(f_k)$ are left dependent over $P^e$.
\end{lm}
\Proof
Let $F=F(z_1,z_2,\ldots,z_k)$ be a nonzero element of $T=k\{z_1,z_2,\ldots,z_k\}$ with the minimal degree such that
$F(f_1,f_2,\ldots,f_k)=0$. Suppose that
\bes
H(F)=u_1 H(z_1)+u_2 H(z_2)+\ldots + u_k H(z_k)
\ees
in $\Omega_T$. We may assume that $u_1=u_1(z_1,z_2,\ldots,z_k)\neq 0$. Note that $\deg\,u_1 < \deg\,F$. Consequently,
\bes
0=H(F(f_1,f_2,\ldots,f_k))= u_1' H(f_1)+u_2' H(f_2)+\ldots + u_k' H(f_k),
\ees
where $u_i'=u_i(f_1,f_2,\ldots,f_k)$ for all $i$. If $u_1'\neq 0$ then the last equation gives a nontrivial dependency of $H(f_1),H(f_2),\ldots,H(f_k)$.

Suppose that $u_1'= 0$. Note that $u_1=t+w$, where $t\in T$ and $w\in \Omega_T$, since $U(T)=T\oplus \Omega_T$.
Obviously, $t(f_1,f_2,\ldots,f_k)\in P$ and it easily follows from (\ref{f1})--(\ref{f2}) that $w(f_1,f_2,\ldots,f_k)\in \Omega_P$. Then, $t(f_1,f_2,\ldots,f_k)=0$ and $w(f_1,f_2,\ldots,f_k)=0$ since
$0=u_1'=t(f_1,f_2,\ldots,f_k)+w(f_1,f_2,\ldots,f_k)\in P\oplus \Omega_P$. If $t\neq 0$ then this contradicts to the minimality of $\deg\,F$ since $\deg\,t\leq \deg\,u_1<\deg\,F$. If $w\neq 0$ then, continuing the same discussions, we get a nontrivial dependency of $H(f_1),H(f_2),\ldots,H(f_k)$ over $P^e$.  $\Box$

\begin{theor}\label{t5} Let $f$ and $g$ be arbitrary elements of the free Poisson algebra \\ $P=k\{x_1,x_2,\ldots,x_n\}$ in the variables $x_1,x_2,\ldots,x_n$ over a field $k$ of characteristic zero. Then the following conditions are equivalent:

(i) $f$ and $g$ are Poisson dependent;

(ii) $H(f)$ and $H(g)$ are left dependent over $P^e$;

(iii) $f$ and $g$ are polynomially dependent, i.e., they are algebraically dependent in the polynomial algebra $P$;

(iv) there exists $a\in P$ such that $f,g\in k[a]$;

(v) $\{f,g\}=0$ in $P$.
\end{theor}
\Proof By Lemma \ref{l12}, (i) implies (ii). The conditions (iii), (iv), and (v) are equivalent \cite{MakarU2,Zaks}.
Obviously, (iii) implies  (i).

To prove the theorem it is sufficient to show that (ii) implies (iii). Note that $f=0$ if and only if $H(f)=0$. Suppose that $(H(f),H(g))\neq 0$ and
\bes
u H(f)+ v H(g)=0, \ \ (u,v)\neq 0.
\ees
Then obviously
\bes
\overline{\overline{u} \overline{H(f)}+ \overline{v} \overline{H(g)}}=0, \ \ (\overline{u},\overline{v})\neq 0,
\ees
or equivalently,
\bes
\varphi(u) \varphi(H(f))+ \varphi(v) \varphi(H(g))=0, \ \ (\varphi(u),\varphi(v))\neq 0
\ees
in the algebra $B=P\otimes_k A$,
where $\varphi : P^e\rightarrow gr\,P^e$ is the gradation mapping (\ref{f9}).
So, $\varphi(H(f))$ and $\varphi(H(g))$ are left dependent over $B$.

Let $l=l(x_1,x_2,\ldots,x_n)$ be an arbitrary element of the free Lie algebra \\ $g=Lie<x_1,x_2,\ldots,x_n>$. Then,
$H_l=h_l=l(h_{x_1},h_{x_2},\ldots,h_{x_n})$ by (\ref{f1}). Hence  $H_l=l(y_1,y_2,\ldots,y_n)\in L$. For every $i\geq 1$ denote by $\partial_{e_i}$ the usual partial derivation of the polynomial algebra $P$ in the variables (\ref{f5}). It can be easily checked that
\bes
H(a)= \sum_{i\geq 1} \partial_{e_i}(a) H(e_i) \in P L
\ees
for every $a\in P$. Consequently, $\varphi(H(f)), \varphi(H(g))\in P\otimes_k L$.

So, the homogeneous nonzero elements $\varphi(H(f))$ and $\varphi(H(g))$ of $P\otimes_k L$ are left dependent over $B$. Then $Deg\,\varphi(H(f))$=$Deg\,\varphi(H(g))$ by Lemma \ref{l11}. Recall that $B$ is a free associative algebra over $P$. Hence there exists nonzero elements $\lambda, \mu \in P$ such that $\lambda \varphi(H(f))=\mu \varphi(H(g))$ or equivalently, $\lambda \overline{H(f)}=\mu \overline{H(g)}$ and $hdeg\,(\lambda H(f)-\mu H(g))< hdeg\,H(f)=hdeg\,H(f)$.

Using Lemma \ref{l7}, it is easy to show that $H(f)$ and $\lambda H(f)-\mu H(g)$ are left dependent over $P^e$ again. If $\lambda H(f)-\mu H(g)\neq 0$ then $\varphi(H(f))$ and $\varphi(\lambda H(f)-\mu H(g))$ are homogeneous nonzero elements of $P\otimes_k L$ left dependent over $B$. This contradicts to the statement of the Lemma \ref{l11} since $Deg\,\varphi(\lambda H(f)-\mu H(g))<Deg\,\varphi(H(f))$. Consequently, $\lambda H(f)-\mu H(g)=0$. Then,
\bes
\sum_{i\geq 1}(\lambda \partial_{e_i}(f)-\mu \partial_{e_i}(g)) H(e_i)=0
\ees
and hence
\bes
\lambda \partial_{e_i}(f)-\mu \partial_{e_i}(g)=0
\ees
for all $i\geq 1$. It is well-known (see, for example \cite{Umi24}) that in this case $f$ and $g$ are algebraically dependent in the polynomial algebra $P$.
$\Box$

For every $p\in P=k\{x_1,x_2,\ldots,x_n\}$ the Fox derivatives  $\frac{\partial p}{\partial x_i}$ (see \cite{Um11,Um}) are uniquely defined by
\bes
H(p)= \frac{\partial p}{\partial x_i} h_{x_1} + \frac{\partial p}{\partial x_i} h_{x_2} + \ldots + \frac{\partial p}{\partial x_i} h_{x_n}, \ \ \ \frac{\partial p}{\partial x_i}\in P^e,
\ees
for all $1\leq i \leq n$. For every endomorphism $\psi$ of the free Poisson algebra $P$ we define the Jacobian matrix $J(\psi)=[u_{ij}]$ with $u_{ij}=\frac{\partial \psi(x_i)}{\partial x_j}$ for all $1\leq i, j \leq n$. It is easy to show \cite{Umi24} that $J(\psi)$ is invertible over $P^e$ if $\psi$ is an automorphism. The reverse statement is an analogue of the classical Jacobian Conjecture for free Poisson algebras.

\begin{theor}\label{t6} Let $\psi$ be an endomorphism of the free Poisson algebra $P=k\{x,y\}$ in the variables $x,y$ over a field $k$ of characteristic $0$. If $J(\psi)$ is invertible over $P^e$ then $\psi(x), \psi(y)\in k[x,y]$.
\end{theor}
\Proof Put $\psi(x)=f$ and $\psi(y)=g$. Note that $J(\psi)$ is invertible over $P^e$ if and only if $\Omega_P$ is the free $P^e$ module with basis $H(f)$ and $H(g)$. It is sufficient to prove that $hdeg\,(H(f))=hdeg\,(H(g))=1$.

Suppose that $hdeg\,H(f)+hdeg\,H(g)\geq 3$. Note that $\Omega_P=P^e H(f)+P^e H(g)=P^e h_x+P^e h_y$. Consequently, $P^e H(f)+P^e H(g)$ contains two elements $h_x$ and $h_y$ of $h$-degree $1$. Hence there exists $(u,v)\neq 0$ such that
\bes
\overline{u \overline{H(f)}+ v \overline{H(g)}}=0.
\ees
As in the proof of Theorem \ref{t5}, $hdeg\,H(f)=hdeg\,H(g)$ and there exist $0\neq \lambda, \mu \in P$ such that $hdeg\,(\lambda H(f)-\mu H(g))< hdeg\,H(f)$. Put $T=\lambda H(f)-\mu H(g)$. Note that $hdeg\,H(f)+hdeg\,T\geq 3$. By Lemma \ref{l8}, it is not difficult to find a nonzero $\eta\in P$ such that
$\eta h_x, \eta h_y\in P^e H(f)+P^e T$. Again, as in the proof of Theorem \ref{t5}, we get $hdeg\,H(f)=hdeg\,T$. This is a contradiction.
$\Box$

Using Jung's Theorem \cite{Jung} and Theorem \ref{t6} we get
\begin{co}\label{c4}
Automorphisms of free Poisson algebras in two variables over a
field of characteristic zero are tame.
\end{co}
This is the third proof of this theorem. The first one was given in \cite{MLTU} by studing the locally nilpotent derivations and the second was given in \cite{MLU} as a corollary of the Freiheitssatz.
\begin{co}\label{c5}
The two-dimensional Jacobian Conjecture for free Poisson algebras is equivalent to the two-dimensional Jacobian Conjecture for polynomial algebras in characteristic zero.
\end{co}

\section{Comments and Problems}

\hspace*{\parindent}

In Section 3 we proved that the universal enveloping algebras of the Poisson symplectic algebra $P_n$ and the Weyl algebra $A_n$ are isomorphic. It is not difficult to show that this result is true also for fields of positive characteristic. A. Belov-Kanel and M. Kontsevich \cite{BKK1} formulated the next problem.
\begin{prob}\label{prob01} The automorphism group of the Weyl algebra of index $n$ is isomorphic to the group of polynomial symplectomorphisms of a $2n$-dimensional affine space, i.e.,
\bes
\Aut\,A_n \simeq \Aut\,P_n.
\ees
\end{prob}
This problem was posed in \cite{BKK1} for fields of characteristic zero but it makes sense in positive characteristic also \cite{AdEssen,BKK2,Tsuchimoto1,Tsuchimoto2}.

Every automorphism $\varphi$ of $P_n$ can be uniquely extended to an automorphism $\varphi^*$ of the universal enveloping algebra $P_n^e$ by $x\mapsto \varphi(x), h_x\mapsto h_{\varphi(x)}$ for all $x\in P_n$. Similarly, every automorphism $\psi$ of $A_n$ can be uniquely extended to an automorphism $\psi^\circ$ of the universal enveloping algebra $A_n^e$ by $x\otimes 1\mapsto \psi(x)\otimes 1, 1\otimes x'\mapsto 1\otimes\psi(x)'$ for all $x\in A_n$. By means of $\varphi\mapsto\varphi^*$ and $\psi\mapsto\psi^\circ$ we can identify the groups of automorphisms
$Aut\,P_n$ and $Aut\,A_n$ with the corresponding subgroups of $Aut\,P_n^e$ and $Aut\,A_n^e$, respectively. By means of the canonical isomorphism $\theta$ from Section 3 we can identify $P_n^e$ and $A_n^e$ and consider $Aut\,P_n$ and $Aut\,A_n$ as subgroups of $Aut\,A_{2n}$.

\begin{prob}\label{prob02} Is it true that $\Aut\,P_n$ and $\Aut\,A_n$ are conjugate?
\end{prob}

First of all it is interesting to know the answer to the question: Are $\Aut\,P_n$ and $\Aut\,A_n$ conjugate in $Aut\,A_{2n}$? It seems the structure of an isomorphism between $\Aut\,P_n$ and $\Aut\,A_n$ is very complicated if it exists.

The structure theory of Poisson algebras is not developed yet. In characteristic zero the universal enveloping algebras of simple Poisson symplectic algebras, the Weyl algebras, are also simple.
\begin{prob}\label{prob03}
Is it true that every Poisson algebra $P$ is (Poisson) simple if and only if its universal enveloping algebra $P^e$ is simple?
\end{prob}
\begin{prob}\label{prob04}
Describe all simple Poisson brackets on the polynomial algebra $k[x,y,z]$ in three variables.
\end{prob}

It is well known (see \cite{Um11,Um}) that the universal enveloping algebras of free Lie algebras are free associative algebras and  every left ideal of a free associative algebra is a free left module \cite{Cohn}. These results are closely related to \cite{Um} classical theorems on free Lie algebras: Subalgebras of free Lie algebras are free \cite{Shir1,Witt}; Automorphisms of finitely generated free Lie algebras are tame \cite{Cohn2}.

We say that the {\em subalgebra membership problem} is decidable for an algebra $A$ if there is an effective procedure that defines for any element $a\in A$  and for any finitely generated subalgebra $B$ of $A$ whether  $a$ belong to $B$ or not. Some relations between the distortion of subalgebras and the subalgebra membership problem were established in \cite{BO} for  polynomial, free associative, and free Lie algebras. It can be easily derived from \cite{Shir1,Witt} that the subalgebra membership problem is decidable for free Lie algebras. Moreover, finitely generated subalgebras of free Lie algebras are residually finite \cite{Um6}. The subalgebra membership problem is still open in the case of free Poisson algebras.
\begin{prob}\label{prob1}
Is the subalgebra membership problem decidable for free Poisson algebras?
\end{prob}

The subalgebra membership problem for free associative algebras is undecidable \cite{UUU,UUU1}. In fact, if $A=k\langle x_1,\ldots,x_n\rangle$ is a free associative algebra, then the structure of the left ideals of the universal enveloping algebra
$A^e=A\otimes_k A^{op}$ is very difficult \cite{UUU}. The left ideal membership problem is algorithmically undecidable and the left dependency of a finite set of elements of $A^e$ is algorithmically unrecognizable \cite{UUU}. The next problem is closely related to Problem \ref{prob1}.
\begin{prob}\label{prob2}
Is the left ideal membership problem decidable over the universal enveloping algebras of free Poisson algebras?
\end{prob}

By Theorem \ref{t4}, the left dependency of a finite set of elements over the universal enveloping algebras of free Poisson algebras is algorithmically recognizable. This is a positive result in the direction of solving Problems \ref{prob1} and \ref{prob2}. These problems are also related to the next problem.
\begin{prob}\label{prob3}
Is the freeness of a finite set of elements of free Poisson algebras algorithmically recognizable?
\end{prob}
Recall that the freeness of a finite set of elements is algorithmically recognizable in the case of free Lie algebras (see, for example \cite{Um9}) and unrecognizable in the case of free associative algebras \cite{UUU}. In Lemma \ref{l13} we proved that the Poisson  dependency implies the dependency of differentials.

\begin{prob}
 Let $f_1,f_2,\ldots,f_k$ be arbitrary elements of the free Poisson algebra
 $P$ over a field $k$ of characteristic $0$. Is it true that the left dependency of
 $H(f_1),H(f_2),\ldots,H(f_k)$ over $P^e$ imply the Poisson dependency of  $f_1,f_2,\ldots,f_k$.
\end{prob}

In relation with Corollary \ref{c5} we formulate one more question.
\begin{prob}\label{prob4}
Is the Jacobian conjecture for free Poisson algebra equivalent to the classical Jacobian conjecture?
\end{prob}
The classical Jacobian conjecture is stably equivalent to homogeneous degree 3 case \cite{Bass,Essen}.
\begin{prob}\label{prob5}
Is the Jacobian conjecture for free Poisson algebras stably reducible to homogeneous degree 3 case?
\end{prob}
An analogue of the Jacobian Conjecture is true for free Lie algebras \cite{Reu,Shpil,Um9,ZM}, for free associative algebras \cite{DL,Sch}, and for free nonassociative algebras \cite{Um11,Yag}.

\bigskip

\begin{center}
{\bf\large Acknowledgments}
\end{center}

\hspace*{\parindent}

I am grateful to Max-Planck Institute f\"ur Mathematik for
hospitality and excellent working conditions, where part of this work has been done. I am also grateful to V. Gorbounov and L. Makar-Limanov  for very helpful discussions and comments.

\end{document}